\documentclass[11pt]{amsart}
\usepackage{amsmath,amssymb,amsthm,amscd}
\usepackage[mathscr]{eucal}
\usepackage[hidelinks,unicode]{hyperref}
\usepackage{xcolor}

\numberwithin{equation}{section}
\allowdisplaybreaks[1]
\renewcommand{\leq}{\leqslant}
\renewcommand{\geq}{\geqslant}
\theoremstyle{plain}

\def\dsum{\displaystyle\sum}

\newtheorem{prop}{Proposition}
\newtheorem{theo}[prop]{Theorem}
\newtheorem{lemm}[prop]{Lemma}
\newtheorem{coro}[prop]{Corollary}
\theoremstyle{definition}
\newtheorem{defi}[prop]{Definition}
\theoremstyle{remark}
\newtheorem{rema}[prop]{Remark}
\newcommand{\R}{\mathbb R}
\newcommand{\Sym}{\operatorname{Sym}}
\newcommand{\tr}{\operatorname{tr}}
\newcommand{\diag}{\operatorname{diag}}

\newcommand{\Om}{\Omega}
\newcommand{\OO}{\overline\Omega}
\newcommand{\II}{\mathrm{II}}
\newcommand{\norm}[1]{\left\|#1\right\|}
\newcommand{\abs}[1]{\left|#1\right|}

\title[The Dirichlet problem for sum Hessian equations]
{The Dirichlet problem for sum Hessian equations on uniformly $(k-1)$-convex domains}

\begin{document}
\author{Changyu Ren}
\address{School of Mathematical Science\\
Jilin University\\ Changchun\\ China}
\email{rency@jlu.edu.cn}
\author{Ziyi Wang}
\address{School of Mathematical Science\\
Jilin University\\ Changchun\\ China}
\email{ziyiw23@mails.jlu.edu.cn}
\thanks{Research of the first author is supported by NSFC Grant No. 12571216.}
\begin{abstract}
In this paper, we mainly study the sum Hessian equation
with general boundary values in the case $2\leq k\leq n$.
Under the conditions that the boundary is uniformly 
$(k-1)$-convex, the right-hand side is monotone non-decreasing in 
$u$, convex in the gradient, and has a suitably small linear gradient growth,
we prove the existence of a unique smooth
admissible solution.
The proof constructs the global subsolution and gradient
barriers by using strictly $k$-convex function, and derives
the boundary second-order estimates from the tangential linearizations and a
duality inequality for the boundary curvature.
\end{abstract}
\maketitle

\section{Introduction}\label{sec:intro}

We consider the following Dirichlet problem:
\begin{equation}\label{eq:problem}
\begin{cases}
S_k(D^2u):=\sigma_k(D^2u)+a\sigma_{k-1}(D^2u)=f(x,u,Du),&x\in\Om,\\
u=\phi,&x\in\partial\Om,
\end{cases}
\end{equation}
where $u$ is an unknown function defined on $\Omega$, $f$ is a
positive function, $2\leq k\leq n$, and $a\geq0$.
Let $\sigma_k(D^2u)=\sigma_k(\lambda(D^2u))$ denote the $k$-th
elementary symmetric function of the eigenvalues of the Hessian
matrix $D^2u$. Namely, for
$\lambda=(\lambda_1,\cdots,\lambda_n)\in\mathbb R^n$,
\begin{equation*}
\sigma_k(\lambda)=\sum_{1\leq i_1<\cdots<i_k\leq
n}\lambda_{i_1}\cdots\lambda_{i_k}.
\end{equation*}
And $\sigma_0=1$.
The admissibility condition required is
\begin{equation}\label{eq:admissible}
\lambda(D^2u(x))\in\Gamma_{k-1}^{n},\quad x\in\OO,
\end{equation}
where
\begin{equation*}
\Gamma_j^N=\{\lambda\in\R^N\big|\sigma_l(\lambda)>0,\ 1\leq l\leq j\},
\end{equation*}
and throughout the paper $\nu$ denotes the outer unit normal.

Hessian equations form a fundamental class of fully nonlinear elliptic
equations, closely related to curvature problems in geometric analysis.
When $a=0$, problem \eqref{eq:problem} reduces to the classical $k$-Hessian
equation.
For the right hand side function $f$ not depending on the gradient term,
Caffarelli-Nirenberg-Spruck \cite{CNS} solved the Dirichlet
problem.
If the right-hand side depends on the gradient, establishing the $C^2$ estimates of the solution
is a long-standing problem. For some works related to
the $C^2$ estimates, see \cite{ChouWang, GuanLiLi, JW, Q, WY}.

In recent years, the fully nonlinear equations derived from linear
combinations of elementary symmetric functions have received increasing
attention. For example, the following special Lagrangian equation
$${\rm Imdet}(\delta_{ij}+{\rm i}u_{ij})=\dsum_{k=0}^{[(n-1)/2]}(-1)^{k}\sigma_{2k+1}
(\lambda(u_{ij}))=0$$ was derived by Harvey and Lawson in their
study of the minimal submanifold problem \cite{HarveyLawson}. The fully
nonlinear partial differential equations
$$P_m(u_{ij})=\sum_{k=0}^{m-1}(l_k^+)^{m-k}(x)\sigma_k(u_{ij})
$$
were studied by Krylov \cite{Krylov} and Dong \cite{Dong2006}. The following
curvature equations and complex Hessian equations
$$\dsum_{s=0}^k\alpha_s\sigma_s(\kappa(X))=f(\chi,\nu(X)),\quad
X\in M,$$
$$\dsum_{s=0}^k\alpha_s\sigma_s(\lambda)=f(z,u,Du)$$
were studied by Li-Ren-Wang \cite{LiRenWang} and Dong \cite{Dws}. For the
following curvature equations
$$\sigma_k(W_u(x))+\alpha(x)
\sigma_{k-1}(W_u(x))=\dsum_{l=0}^{k-2}\alpha_l(x)\sigma_l(W_u(x)),\quad
x\in\mathbb{S}^n,$$ Guan-Zhang \cite{GuanZhang} and Zhou \cite{Zhou}
studied the curvature estimates and the interior gradient estimates
respectively.

For two-term sum equations over the standard Hessian, Liu-Ren
\cite{LiuRen} studied Pogorelov type second-order estimates for the zero
Dirichlet boundary value problem.
In recent years, the related estimates have been further developed for weaker
classes of solutions.
Liang-Yan-Zhu \cite{LiangYanZhuPogorelov} established Pogorelov type
estimates under $(k-1)$-convexity and semiconvexity conditions.
Jin-Xu \cite{JinXu} also studied Pogorelov type estimates for admissible
solutions with semi-convexity when the right-hand side depends on the gradient.
For combinations involving elementary symmetric functions of several
different orders, Dong-Xu-Zhang \cite{DongXuZhang2026} studied Pogorelov type interior estimates under a
real-root condition and the corresponding admissibility assumptions, and
discussed the dynamic semiconvexity in the cases $2\leq k\leq n-2$.

Concerning boundary problems, Jia-Yang \cite{JiaYang} established
Pogorelov type estimates for $k$-convex solutions up to flat boundaries for
sum Hessian equations with right-hand side $f(x,u)$ in half space.
Liang-Yan-Zhu \cite{LiangYanZhuNeumann} established
second-order a priori estimates and existence results for the Neumann problem.
These works involve different boundary conditions, different dependencies of
the right-hand side and different admissibility classes.
In order to treat general Dirichlet boundary values for problem
\eqref{eq:problem}, it remains necessary to establish gradient bounds,
boundary second-order bounds, and uniform estimates along the path of
continuity that are adapted to its structure.

Regarding $S_k$, when $a>0$, the operator $S_k^\frac{1}{k}$ with respect to $D^2u$ is no longer homogeneous of degree one,
and generally allows for $\sigma_k(D^2u)\leq0$. This makes the treatment of $S_k$ considerably more challenging.
Guan and Jiao~\cite{GuanJiao2015} studied the Dirichlet problem
\begin{equation*}
\begin{cases}
f\bigl(\lambda(\nabla^2 U+A(X,U,\nabla U))\bigr)
  = \psi(X,U,\nabla U), & X\in M,\\
U=\varphi, & X\in\partial M,
\end{cases}
\end{equation*}
on a Riemannian manifold $M$ with smooth boundary.
Our two-term sum operator fits this general framework on $M=\Omega\times\mathbb S^1$ by adjoining a fixed eigenvalue and denoting
$S_k(D^2u)$ as $\sigma_k(\lambda(D^2u),a)$, but
the assumptions on the right-hand side differ.
Ma-Wang-Zhu~\cite{MaWangZhu} also use the identity
$\sigma_k(\lambda,a)=\sigma_k(\lambda)+a\sigma_{k-1}(\lambda)$
to establish ellipticity and concavity for a sum operator acting on
the first Newton transform of the second fundamental form of a
closed star-shaped hypersurface. We also mention that Guan~\cite{Guan1994} used admissible strict subsolutions to solve
a class of fully nonlinear Dirichlet problems with gradient-dependent
right-hand sides. We also use an admissible strict subsolutions in our proof.

We now present the conditions on the domain and the main results.

\begin{defi}\label{def:geometry}
Call $\partial\Om$ uniformly $(k-1)$-convex if there exists
$\delta>0$ such that
\begin{equation}\label{eq:geometry}
\lambda(\II(x)-\delta I)\in\Gamma_{k-1}^{n-1},\quad x\in\partial\Om,
\end{equation}
\end{defi}
where the second fundamental form is taken as $\II(X,Y)=\langle D_X\nu,Y\rangle$ so that the
principal curvatures of a sphere are positive. Let $\kappa_1(x),\cdots,\kappa_{n-1}(x)$ denote principal curvatures of the boundary and $\kappa(x)=(\kappa_1(x),\cdots,\kappa_{n-1}(x)):=\lambda(\II(x))$.
Since the boundary is compact, this is equivalent to the range of the
principal curvature vector being compactly contained in $\Gamma_{k-1}^{n-1}$.
When $k=n$, it is precisely uniform convexity.

In \cite{MaZhang}, there is a defining function $\rho\in C^\infty(\OO)$
satisfying
\begin{equation*}
\rho<0\ \text{in }\Om,\quad
\rho=0,\quad \rho_\nu=1\ \text{on }\partial\Om,\quad
D^2\rho\in\Gamma_k^n\ \text{in }\OO.
\end{equation*}
Caffarelli-Nirenberg-Spruck \cite{CNS} proved the existence of this smoothly and strictly $k$-convex defining function $\rho$ on $\bar{\Omega}$.
Denote
\begin{equation}\label{eq:b0}
\beta_0=\frac{\min\limits_{\bar{\Omega}}\sigma_k^\frac{1}{k}(D^2\rho)}{\norm{D\rho}_\infty}>0,
\end{equation}
which will be used in the following theorem.
\begin{theo}\label{thm:main}
Let $\Om\subset\R^n$ be a bounded smooth domain whose boundary satisfies
\eqref{eq:geometry}.
Let $a\geq0$ be a constant, $\phi\in C^\infty(\partial\Om)$, and
$f\in C^\infty(\OO\times\R\times\R^n)$ with $f>0$. Denote
\begin{equation*}
g(x,z,p)=f^\frac{1}{k}(x,z,p).
\end{equation*}
Assume that for $x\in\OO$, $z\leq \max_{\partial\Om}\phi$ and $p\in\R^n$:
\begin{enumerate}
\renewcommand{\labelenumi}{\textup{(H\arabic{enumi})}}
\item $g$ is a convex function with respect to $p$ and $g_z\geq0$;
\item there exist constants $C_0>0$ and $0\leq\beta<\beta_0$ such that
\begin{equation}\label{eq:growth}
g(x,z,p)\leq C_0+\beta\abs{p},
\end{equation}
\item for every compact interval $J\subset(-\infty,\max_{\partial\Om}\phi]$,
\begin{equation}\label{eq:gx}
\sup_{\OO\times J\times\R^n}\abs{D_xg}<\infty.
\end{equation}
\end{enumerate}
Then problem \eqref{eq:problem} exists a unique solution
$u\in C^\infty(\OO)$ satisfying \eqref{eq:admissible}.
\end{theo}
\begin{rema}\label{rem:gx}
Condition \eqref{eq:gx} is used only before gradient bounds are available.
Once the $C^1$ bound is established, the other derivatives of finite order
need only be bounded on a fixed compact set.
\end{rema}
\begin{rema}\label{rem:gradient-convexity}
In the proof of Theorem~\ref{thm:main}, the assumption in
\textup{(H1)} that $g=f^\frac{1}{k}$ is convex in $p$ is used
only to derive the global $C^2$ estimates from
the boundary $C^2$ estimates.
\end{rema}

Without assuming the
convexity of $g$ in $p$, there have been several results for the
classical $k$-Hessian equation.
Guan-Ren-Wang~\cite{GuanRenWang} first obtained the global $C^2$ estimates for $k=2$.
Ren-Wang~\cite{RenWang1, RenWang2} then dealt with the problem for $k=n-1, n-2$.
Subsequently, Yan~\cite{Yan2026Curvature} got the global $C^2$ estimates when $n\geq 3$ and $n/2\leq k<n$.
The existence of solutions can still be obtained in the above works.

For sum Hessian equation without the convexity of $g$ in $p$, when $k=n$, the Pogorelov type estimate has been established by Li-Ren \cite{LiRen}. While $k=2$ and $k=3$,
Liu-Ren \cite{LiuRen} have obtained the Pogorelov estimates. Besides, for general $k$, Liang-Yan-Zhu \cite{LiangYanZhuPogorelov} and Jin-Xu \cite{JinXu} have proved Pogorelov
estimates under the condition that the solution is semi-convex.

\begin{coro}\label{cor:no-gradient}
If $f=f(x,z)>0$ is smooth with $f_z\geq0$, then on a uniformly
$(k-1)$-convex domain as above, for every
$\phi\in C^\infty(\partial\Om)$, problem \eqref{eq:problem} exists a unique
smooth $\Gamma_{k-1}$-admissible solution.
\end{coro}
\begin{proof}
In this case $g_p=0$, $D_{pp}^2g=0$, and
\begin{equation*}
g(x,z)\leq g(x,\max_{\partial\Om}\phi)\leq\max_{\OO}g(x,\max_{\partial\Om}\phi),\quad z\leq \max_{\partial\Om}\phi.
\end{equation*}
Hence we take $\beta=0$ and condition \eqref{eq:gx} holds by smoothness and
compactness.
\end{proof}

The organization of our paper is as follows.
In Section 2, we give the preliminary knowledge.
The $C^0$ estimates are given in Section 3.
In Section 4, we give the gradient estimates.
The boundary tangential and mixed $C^2$ estimates
is contained in Section 5.
In Section 6, we give the boundary normal $C^2$ estimates
and briefly introduce how to obtain the global $C^2$
estimates.
In Section 7 we used the method of continuity to complete the proof of Theorem \ref{thm:main}.

\section{Preliminary}
We recall some basic definition and properties as follows.
\begin{prop}\label{lem2.1}We have
\par
(1) $S_k(\lambda)=\lambda_i S_{k-1}(\lambda|i)+S_k(\lambda|i),\quad
i=1,2,\cdots,n$;
\par
(2) $S_k^{pp}(\lambda):=\frac{\partial
S_k(\lambda)}{\partial\lambda_p}=\sigma_{k-1}(\lambda|p)+a\sigma_{k-2}(\lambda|p)=S_{k-1}(\lambda|p),\quad
p=1,2,\cdots,n$;
\par
(3) $\displaystyle\sum_{i=1}^n
S_k(\lambda|i)=(n-k)S_k(\lambda)+a\sigma_{k-1}(\lambda)$;
\par
(4) $\displaystyle\sum_{i=1}^n\lambda_i
S_{k-1}(\lambda|i)=kS_k(\lambda)-a\sigma_{k-1}(\lambda)$.
\par
(5) (Newton-Maclaurin) If~$\lambda=(\lambda_1,\cdots,\lambda_n)\in \Gamma_k^n$, then
\begin{equation*}
\frac{\sigma_1(\lambda)}{C_n^1}\geq(\frac{\sigma_2(\lambda)}{C_n^2})^\frac{1}{2}\geq\cdots\geq(\frac{\sigma_k(\lambda)}{C_n^k})^\frac{1}{k};
\end{equation*}
\end{prop}
\begin{proof}
(1)-(4) can be found in \cite{LiuRen}.
\end{proof}
\par
For $1\leq m\leq N$, $\Gamma_m^N$ is an open convex cone, invariant under the
addition of positive semidefinite matrices.

\par
The following lemma comes from \cite{CNS}.
\begin{lemm}\label{BC}
For $B,C\in\Gamma_m^N$,
\begin{equation}\label{eq:support}
\Phi_m^{ij}(B)C_{ij}\geq\Phi_m(C)>0,
\end{equation}
where $\Phi_m(B)=\sigma_m(B)^{1/m}$.
\end{lemm}
Caffarelli-Nirenberg-Spruck \cite{CNS} have also proved the concavity and 1-homogeneity of $\sigma_m(B)^{1/m}$.
\par
Denote
\begin{equation*}
F(A):=S_k^\frac{1}{k}(\lambda(A)),\quad F^{ij}:=\frac{\partial F}{\partial A_{ij}},\quad F^{ij,rs}:=\frac{\partial^2 F}{\partial A_{ij}\partial A_{rs}},\quad F_i:=\frac{\partial F}{\partial \lambda_i}
\end{equation*}
and define the positive admissible branch
\begin{equation*}
\widetilde\Gamma_k^n=\{A\in\Sym(n)\big|\lambda(A)\in\Gamma_{k-1}^n,\ S_k(A)>0\}.
\end{equation*}
Li-Ren-Wang \cite{LiRenWang} have proved that~$F$~is concave in the set
~$\widetilde\Gamma_k^n$.
\begin{lemm}\label{lem:augment}
For $a\geq0$,
\begin{equation*}
A\in\widetilde\Gamma_k^n\quad\Longleftrightarrow\quad\diag(A,a)\in\Gamma_k^{n+1}.
\end{equation*}
Consequently $\widetilde\Gamma_k^n$ is an open convex set.
\end{lemm}
\begin{proof}
We have
\begin{equation*}
\sigma_j(\lambda(A),a)=\sigma_j(A)+a\sigma_{j-1}(A),\quad1\leq j\leq k.
\end{equation*}
If $A\in\widetilde\Gamma_k^n$, then the right-hand side is positive for $j\leq k-1$ and equals
$S_k(A)>0$ for $j=k$. Hence the augmented matrix lies in $\Gamma_k^{n+1}$.
Conversely, deleting the last row and column of the augmented matrix yields
$A\in\Gamma_{k-1}^n$, and its $k$-th elementary symmetric function is
$\sigma_k(\diag(A,a))=S_k(A)>0$.
\end{proof}

\begin{rema}\label{rem:nonhomogeneous}
When $a>0$ is fixed, $\widetilde\Gamma_k^n$ is not a cone in general, and $F$ is not
homogeneous of degree one in $A$. In the sequel all homogeneity
computations are carried out on the $(n+1)$-dimensional augmented operator.
\end{rema}

The following lemma is a direct consequence of Lemma \ref{BC}.

\begin{lemm}\label{lem:linear}
If $A\in\widetilde\Gamma_k^n$ and $B\in\Gamma_k^n$, then
\begin{equation}\label{eq:linear}
F^{ij}(A)B_{ij}\geq\sigma_k^\frac{1}{k}(B).
\end{equation}
In particular,
\begin{equation}\label{eq:trace-rho}
\tr F:=\sum_iF^{ii}\geq\binom nk^\frac{1}{k}.
\end{equation}
\end{lemm}

\begin{lemm}\label{lem:largest}
Let $\lambda_1$ be the largest eigenvalue of $A\in\widetilde\Gamma_k^n$, and suppose
$\lambda_1\geq a$. Then
\begin{equation}\label{eq:largest}
\lambda_1F_1\geq\frac{F}{n+1},\quad
\sum_iF_i\lambda_i^2\geq\frac{F\lambda_1}{n+1}.
\end{equation}
\end{lemm}
\begin{proof}
Suppose $\lambda_1,\ldots,\lambda_n$ are eigenvalues of $A$, $\mu=(\lambda_1,\ldots,\lambda_n,a)$.
Then 
\begin{equation*}
S_k(A)=\lambda_1\sigma_{k-1}(\mu|1)+\sigma_k(\mu|1)
\end{equation*}
and
\begin{equation*}
F_1=\frac{1}{k}S_k^{\frac{1}{k}-1}\sigma_{k-1}(\mu|1).
\end{equation*}
If $\sigma_k(\mu|1)\leq0$, then $S_k\leq\lambda_1\sigma_{k-1}(\mu|1)$, thus
\begin{equation*}
\lambda_1F_1\geq \frac{F}{k}\geq \frac{F}{n+1}.
\end{equation*}
If $\sigma_k(\mu|1)>0$, then $\mu|1\in\Gamma_k^{n}$, so that $\sigma_{k-1}(\mu|1j)>0$.
Since $\mu_j\leq\lambda_1$ for all $j\ne1$,
\begin{equation*}
k\sigma_k(\mu|1)=\sum_{j\ne1}\mu_j\sigma_{k-1}(\mu|1j)
\leq\lambda_1\sum_{j\ne1}\sigma_{k-1}(\mu|1j)
=(n-k+1)\lambda_1\sigma_{k-1}(\mu|1).
\end{equation*}
Therefore,
\begin{equation*}
S_k\leq\frac{n+1}{k}\lambda_1\sigma_{k-1}(\mu|1),
\end{equation*}
that is,
\begin{equation*}
\lambda_1F_1\geq \frac{F}{n+1}.
\end{equation*}
Moreover, 
\begin{equation*}
\sum_iF_i\lambda_i^2\geq F_1\lambda_1^2\geq\frac{F\lambda_1}{n+1}.
\end{equation*}
Thus we have (\ref{eq:largest}).
\end{proof}

\begin{lemm}\label{lem:compact}
Given $K,s_0>0$, the set
$\{A\in\widetilde\Gamma_k^n\big|\abs A\leq K,\ S_k(A)\geq s_0\}$
is contained in a compact subset of $\widetilde\Gamma_k^n$.
Consequently there exist $0<\theta_0\leq\theta_1<\infty$ such that, on this set,
\begin{equation*}
\theta_0 I\leq(F^{ij}(A))\leq\theta_1 I.
\end{equation*}
\end{lemm}
\begin{proof}
In \cite{LiangYanZhuNeumann}, Liang-Yan-Zhu have proved that the operator $F$ is uniformly elliptic
in this case. The lemma is a direct consequence of this result.
\end{proof}
\section{$C^0$ estimates}\label{sec:path}

In this section, we will construct a subsolution to derive the $C^0$ estimates.

\par
Let $h\in C^\infty(\OO)$ be the harmonic
extension of $\phi$ satisfying
\begin{equation*}
\begin{cases}
\Delta h=0\ \text{in }\Om,\\
h=\phi\ \text{on }\partial\Om.
\end{cases}
\end{equation*}
The maximum principle gives $h\leq \max\limits_{\partial\Omega}\phi$. For $A>0$, define
\begin{equation*}
\underline{u}=A\rho+h.
\end{equation*}
Since $\rho\leq0$, we have $\underline{u}\leq h\leq \max\limits_{\partial\Omega}\phi$ and $\underline{u}=\phi\ \text{on }\partial\Om$.
Note that the range of $D^2\rho$ is compactly contained in $\Gamma_k^n$. When $A$ is sufficiently large,
\begin{equation*}
D^2\underline{u}=A(D^2\rho+A^{-1}D^2h)\in\Gamma_k^n\subset\widetilde\Gamma_k^n.
\end{equation*}
and
\begin{align*}
\frac{1}{A}F(D^2\underline{u})=&\big[\sigma_k(D^2\rho+A^{-1}D^2h)\\
&+aA^{-1}\sigma_{k-1}(D^2\rho+A^{-1}D^2h)\big]^\frac{1}{k}
\longrightarrow \sigma_k^\frac{1}{k}(D^2\rho).
\end{align*}
\begin{equation*}
F(D^2\underline{u})\geq A(\min\limits_{\bar{\Omega}}\sigma_k^\frac{1}{k}(D^2\rho)-\eta_0),
\end{equation*}
where $\eta_0$ is a positive constant which will be determined later.
From (\ref{eq:growth}), 
\begin{equation*}
g(x,\underline{u},D\underline{u})\leq C_0+\beta\norm{Dh}_\infty+\beta A\norm{D\rho}_\infty.
\end{equation*}
Thus, 
\begin{equation*}
F(D^2\underline{u})-g(x,\underline{u},D\underline{u})\geq A(\min\limits_{\bar{\Omega}}\sigma_k^\frac{1}{k}(D^2\rho)-\eta_0)-C_0-\beta\norm{Dh}_\infty-\beta A\norm{D\rho}_\infty.
\end{equation*}
From (\ref{eq:b0}), choosing
\begin{equation*}
\eta_0=\frac{1}{2}(\min\limits_{\bar{\Omega}}\sigma_k^\frac{1}{k}(D^2\rho)-\beta\norm{D\rho}_\infty)>0
\end{equation*}
and $A$ sufficiently large, we have $F(D^2\underline{u})-g(x,\underline{u},D\underline{u})\geq1$. 
\par
Due to the convexity of $\widetilde\Gamma_k^n$, the ellipticity of $F$, and $g_z\geq0$,
the comparison principle holds. Since $\Delta u>0$,
applying the maximum principle to $u-h$ gives $u\leq h\leq \max\limits_{\partial\Omega}\phi$.
Moreover, $F(D^2\underline{u})-g(x,\underline{u},D\underline{u})\geq0$, so the comparison
principle gives $\underline{u}\leq u$. Therefore we finally obtain
\begin{equation}\label{eq:C0-order}
{\norm u}_{C^0(\bar{\Omega})}\leq C,
\end{equation}
where the constant $C$ depends only on $n$, $k$, $a$, the geometry of $\Omega$, $|\phi|_{C^{2,\alpha}(\partial\Omega)}$, $\rho$, $C_0$ and $\beta$.

\section{$C^1$ estimates}\label{sec:gradient}

We now prove the global gradient bound using a two-point function.
For $D^2u\in\widetilde\Gamma_k^n$, $\lambda=\lambda(D^2u)$, let $\mu=(\lambda,a)\in\Gamma_k^{n+1}$. If $\lambda_i<0$, then
\begin{equation*}
\sigma_{k-1}(\mu|i)
=\sigma_{k-1}(\mu)-\lambda_i\sigma_{k-2}(\mu|i)
\geq\sigma_{k-1}(\mu).
\end{equation*}
where we used $\sigma_{k-2}(\mu|i)>0$.
Combined with Proposition \ref{lem2.1} (5), we obtain
\begin{equation}\label{eq:Fi}
F_i
=\frac{\sigma_{k-1}(\mu|i)}{k\,\sigma_k^{1-\frac{1}{k}}(\mu)}
\geq\frac{\sigma_{k-1}(\mu)}{k\,\sigma_k^{1-\frac{1}{k}}(\mu)}
\geq c(n,k)>0
\quad(\lambda_i<0).
\end{equation}

Set $R\geq\operatorname{diam}(\Omega)$ and 
\begin{equation*}
\omega(r)=K\left(r-\frac{r^\frac{3}{2}}{3\sqrt R}\right),
\quad
P(x,y)=u(x)-u(y)-\omega(|x-y|),
\end{equation*}
where the constant $K>0$ will be determined later.
For $0\leq r\leq R$ we have $Kr/2\leq\omega(r)\leq Kr$.
Since $\underline{u}\leq u\leq h$ and on the
boundary $\underline{u}=u=h$, if $y\in\partial\Omega$, then
\begin{equation*}
\underline{u}(x)-\underline{u}(y)
\leq u(x)-u(y)\leq h(x)-h(y).
\end{equation*}
The functions $h$ and $\underline{u}$ are fixed smooth functions with
Euclidean Lipschitz bounds on $\overline\Omega$.
Therefore, if either $x\in\partial\Omega$ or $y\in\partial\Omega$, then
$|u(x)-u(y)|\leq C|x-y|$ and $P$ is nonpositive
on such pairs of points by choosing $K$ sufficiently large.
Moreover, $P(x,x)=0$. Hence, if the maximum of $P$ on
$\overline\Omega\times\overline\Omega$ is positive, then its maximum point
consists of two distinct interior points $x,y$ and $u(x)>u(y)$.

By (\ref{eq:C0-order}) and $P(x,y)>0$, we have
$|x-y|<4\|u\|_{C^0(\bar{\Omega})}/K$. Thus enlarging $K$ makes $|x-y|$ smaller than a fixed
geometric radius. On this scale, $x,y$ can be joined within $\overline\Omega$ by a curve of
length at most $C|x-y|$. Hence, by (\ref{eq:gx}), for $z\in [\min_{\bar\Omega}\underline{u},\max_{\partial\Om}\phi]$ and $p\in\mathbb R^n$,
\begin{equation*}
|g(y,z,p)-g(x,z,p)|\leq C|x-y|,
\end{equation*}
where $C$ is independent of $p$.

At the maximum point, $D_xP(x,y)=D_yP(x,y)=0$ gives
\begin{equation*}
Du(x)=Du(y)=\omega'(|x-y|)\frac{x-y}{|x-y|}.
\end{equation*}
Using $g_z\geq0$ we get
\begin{align*}
F(D^2u(y))-F(D^2u(x))
&=g(y,u(y),p)-g(x,u(x),p)\\
&\leq g(y,u(y),p)-g(x,u(y),p)
\leq C|x-y|.
\end{align*}
On the other hand, the second-order conditions give, for any $\xi,\eta\in\mathbb R^n$,
\begin{equation*}
\frac{d^2}{dt^2}P(x+t\xi,y+t\eta)\Big|_{t=0}=\xi^{\mathsf T}D^2u(x)\xi-\eta^{\mathsf T}D^2u(y)\eta-\frac{d^2}{dt^2}\omega(|x+t\xi-y-t\eta|)\Big|_{t=0}\leq0.
\end{equation*}
Thus
\begin{align*}
\xi^{\mathsf T}D^2u(x)\xi-\eta^{\mathsf T}D^2u(y)\eta
\leq&\frac{\omega'(|x-y|)}{|x-y|}
\big(|\xi-\eta|^2-[\frac{x-y}{|x-y|}\cdot(\xi-\eta)]^2\big)\\
&+\omega''(|x-y|)[\frac{x-y}{|x-y|}\cdot(\xi-\eta)]^2\\
\leq&\frac{K}{|x-y|}|\xi-\eta|^2.
\end{align*}
In the last inequality we used $\omega''<0$. By taking $\xi=\frac{x-y}{|x-y|}$ and $\eta=0$, the first term in the first inequality equals to 0, and we see that the smallest
eigenvalue $\underline\lambda$ of $D^2u(x)$ satisfies
\begin{equation}\label{eq:lambda*}
\underline\lambda\leq \frac{(x-y)^{\mathsf T}}{|x-y|}D^2u(x)\frac{(x-y)}{|x-y|}\leq\omega''(|x-y|)
=-\frac{K}{4\sqrt{R|x-y|}}<0.
\end{equation}
Take the corresponding unit eigenvector $v$.
By choosing $\xi=\eta+s(\eta\cdot v)v$ with $s=\underline\lambda/(\frac{K}{|x-y|}-\underline\lambda)<0$,
since $D^2u(x)v=\underline\lambda v$, we have
\begin{align*}
\eta^{\mathsf T}(D^2u(y)-D^2u(x))\eta
&\geq \xi^{\mathsf T}D^2u(x)\xi-\eta^{\mathsf T}D^2u(x)\eta
      -\frac{K}{|x-y|}|\xi-\eta|^2\\
&=\bigl[2s\underline\lambda+(\underline\lambda-\frac{K}{|x-y|})s^2\bigr](\eta\cdot v)^2\\
&=s\underline\lambda(\eta\cdot v)^2 .
\end{align*}
Thus $D^2u(y)\geq D^2u(x)+s\underline\lambda v\otimes v$. Note that
\begin{equation}\label{eq:slambda}
-\underline\lambda>s\underline\lambda\geq\frac{K}{16R\bigl(1+\tfrac14\sqrt\frac{r}{R}\bigr)}\geq\frac{K}{20R}.
\end{equation}
Using (\ref{eq:Fi}), (\ref{eq:lambda*}) and (\ref{eq:slambda}) we have
\begin{align*}
F(D^2u(y))-F(D^2u(x))
\geq& F(D^2u(x)+s\underline\lambda v\otimes v)-F(D^2u(x))\\
\geq& c(n,k)s\underline\lambda\\
\geq&\frac{c(n,k)K}{20R}.
\end{align*}
For $K$ sufficiently large, this contradicts $F(D^2u(y))-F(D^2u(x))\leq C|x-y|$ and $|x-y|\leq R$.
Therefore $P(x,y)\leq0$. Interchanging $x,y$, we get
$|u(x)-u(y)|\leq K|x-y|$, and consequently
$\|Du\|_{L^\infty(\Omega)}\leq K$.
Combined with the $C^0$ estimates, we obtain
\begin{equation}\label{eq:C1}
\|u\|_{C^1(\overline\Omega)}\leq C.
\end{equation}
Compared with the constant of the $C^0$ estimates, the above constant also depends on the bounds given by
(H3).

\section{Boundary tangential and mixed $C^2$ estimates}\label{sec:mixed}

We first get a global barrier which will be used in $C^2$ estimates.
\par
By \eqref{eq:C1}, $(x,u,Du)$ lies in a fixed compact set, so the smoothness
of $g$ gives
\begin{equation}\label{eq:gp}
|g_p(x,u,Du)|\leq C.
\end{equation}
\par
Define the linearized operator
\begin{equation*}
L=F^{ij}\partial_{ij}-g_{p_r}\partial_r.
\end{equation*}
Differentiating the equation along the fixed Cartesian direction $x_l$, we
obtain
\begin{equation*}
F^{ij}u_{ijl}=g_{x_l}+g_z u_l
+g_{p_r}u_{rl},
\end{equation*}
therefore
\begin{equation}\label{eq:first-derivative}
(L-g_z)u_l=g_{x_l}.
\end{equation}
\begin{lemm}\label{lem:barrier}
There exist constants $N>0$ and $c_1>0$, such that
for $\widetilde{\rho}=\frac{1-e^{N\rho}}{N}$,
\begin{equation}\label{eq:barrier}
L(\widetilde{\rho})\leq-c_1(1+\tr F)\quad \text{in}~\Omega.
\end{equation}
Moreover, $\widetilde{\rho}\geq0$ on $\overline\Omega$, and $\widetilde{\rho}=0$, $\widetilde{\rho}_\nu=-1$ on
$\partial\Omega$.
\end{lemm}
\begin{proof}
From $\widetilde{\rho}_i=-e^{N\rho}\rho_i$ and
$\widetilde{\rho}_{ij}=-e^{N\rho}(\rho_{ij}+N\rho_i\rho_j)$, we obtain
\begin{equation}\label{Lp}
L\widetilde{\rho}=-e^{N\rho}\bigl[F^{ij}(\rho_{ij}+N\rho_i\rho_j)-g_p\cdot D\rho\bigr].
\end{equation}
Since the range of $D^2\rho$ is compactly contained in $\Gamma_k^n$, we
choose a $\varepsilon>0$ such that
$D^2\rho-2\varepsilon I\in\Gamma_k^n$ on $\overline\Omega$.
Split the matrix as
\begin{equation*}
D^2\rho+N D\rho\otimes D\rho
=(D^2\rho-2\varepsilon I)+\varepsilon I
 +(\varepsilon I+N D\rho\otimes D\rho).
\end{equation*}
The first and third terms both lie in $\Gamma_k^n$.
Using (\ref{eq:linear}) to these two terms
separately, we obtain
\begin{equation}\label{Fp}
F^{ij}(\rho_{ij}+N\rho_i\rho_j)
\geq \min_{\overline\Omega}\sigma_k^\frac{1}{k}(D^2\rho-2\varepsilon I)+\varepsilon\operatorname{tr}F
+\sigma_k^\frac{1}{k}(\varepsilon I+N D\rho\otimes D\rho).
\end{equation}
The eigenvalues of the matrix $\varepsilon I+N D\rho\otimes D\rho$ are
$\varepsilon+N|D\rho|^2$ and $n-1$ copies of $\varepsilon$. Hence
\begin{equation*}
\sigma_k(\varepsilon I+N D\rho\otimes D\rho)
=\binom nk\varepsilon^k
 +\binom{n-1}{k-1}\varepsilon^{k-1}N|D\rho|^2.
\end{equation*}
Therefore
\begin{equation}\label{sigmap}
\sigma_k^\frac{1}{k}(\varepsilon I+N D\rho\otimes D\rho)
\geq
\left[\binom{n-1}{k-1}\varepsilon^{k-1}\right]^\frac{1}{k}
N^\frac{1}{k}\|D\rho\|_\infty^{\frac{2}{k}-1}|D\rho|.
\end{equation}
From (\ref{eq:gp}), (\ref{Lp}), (\ref{Fp}) and (\ref{sigmap}), choosing $N$ sufficiently large, we get
\begin{equation*}
L\widetilde{\rho}\leq-e^{N\rho}(\min_{\overline\Omega}\sigma_k^\frac{1}{k}(D^2\rho-2\varepsilon I)+\varepsilon\operatorname{tr}F)
\leq-c_1(1+\operatorname{tr}F),
\end{equation*}
where
\begin{equation*}
c_1=e^{-N\|\rho\|_\infty}\min\{\min_{\overline\Omega}\sigma_k^\frac{1}{k}(D^2\rho-2\varepsilon I),\varepsilon\}>0.
\end{equation*}
Since $\rho\leq0$, we have $\widetilde{\rho}\geq0$. On the boundary, $\rho=0$ and
$\rho_\nu=1$ give $\widetilde{\rho}=0$ and $\widetilde{\rho}_\nu=-1$, respectively.
\end{proof}
The above inequality holds throughout the entire domain and thus also
applies to any boundary neighborhood.

For tangential vector fields $X,Y$ on the boundary, using $h$ defined in
Section 3, by Gauss formula, we have
\begin{equation*}
u-h=0 \quad\text{on}~\partial\Omega,
\end{equation*}
\begin{equation*}
D^2(u-h)(X,Y)=(u_\nu-h_\nu)\II(X,Y).
\end{equation*}
Thus,
\begin{equation}\label{eq:tangential}
D^2u(X,Y)=D^2h(X,Y)+(u_\nu-h_\nu)\II(X,Y).
\end{equation}
Combined with the boundary gradient bound, we obtain
\begin{equation}\label{eq:tangential-bound}
\abs{D^2u(X,Y)}\leq C\abs X\abs Y\quad\text{on }\partial\Om.
\end{equation}

To estimate mixed derivatives, we use the method in \cite{Guan1994}. Choosing $x_0\in\partial\Om$, translating and rotating the coordinates
so that $x_0=0$ and $\nu(0)=e_n$. Denote $x=(x',x_n)$, $B_r(0):=\{x\in\mathbb{R}^n\big||x|<r\}$ and take the tangential coordinates to be
principal directions. The local boundary can be written as
\begin{equation*}
x_n=\psi_*(x')=-\frac12\sum_{\alpha<n}\kappa_\alpha x_\alpha^2+O(\abs{x'}^3)
\quad \text{on}~\partial\Omega\cap B_r(0).
\end{equation*}
where $\kappa_\alpha$ denotes the principal curvature at $x_0$. Let $w=u-h$, then $w=0$ on the boundary.
For $\alpha<n$, define
\begin{equation}\label{eq:T}
T_\alpha=\partial_\alpha-\kappa_\alpha
(x_\alpha\partial_n-x_n\partial_\alpha)
=(1+\kappa_\alpha x_n)\partial_\alpha-\kappa_\alpha x_\alpha\partial_n
=T_\alpha^l\partial_l.
\end{equation}
where $T_\alpha^l$ are the coefficient components of the vector field.
For the normal vector $\nu$, near $x_0$ we have
\begin{equation*}
\nu=\frac{(-D\psi_*,1)}{\sqrt{1+|D\psi_*|^2}},\quad \nu_\alpha=\kappa_\alpha x_\alpha+O(|x'|^2), \quad \nu_n=1+O(|x'|^2).
\end{equation*}
Thus,
\begin{equation*}
T_\alpha\cdot\nu=\kappa_\alpha x_\alpha-\kappa_\alpha(x_\alpha-x_n\kappa_\alpha x_\alpha)+O(|x'|^2)=O(|x'|^2).
\end{equation*}
Since on the boundary $Dw=w_\nu\nu$,
\begin{equation}\label{eq:Tw-boundary}
\abs{T_\alpha w}= \abs{T_\alpha^l w_l}=\abs{T_\alpha^l w_\nu\nu_l}\leq|u_\nu-h_\nu||T_\alpha\cdot\nu|\leq C\abs x^2\quad
\text{on }\partial\Om\cap B_r(0).
\end{equation}

Differentiating term by term gives
\begin{equation}\label{eq:commutator}
\begin{split}
L(T_\alpha u)={}&F^{ij}T_\alpha^l u_{ijl}+2F^{ij}(\partial_iT_\alpha^l)u_{l j}+F^{ij}(\partial_{ij}T_\alpha^l)u_l\\
&-g_{p_j}(\partial_jT_\alpha^l)u_l-g_{p_j}T_\alpha^l u_{jl}\\
=&T_\alpha^l Lu_l+2F^{ij}(\partial_iT_\alpha^l)u_{l j}
+F^{ij}(\partial_{ij}T_\alpha^l)u_l-g_{p_j}(\partial_jT_\alpha^l)u_l.
\end{split}
\end{equation}
From (\ref{eq:T}) we get
\begin{equation}\label{T}
\partial_{ij}T_\alpha^l=0, \quad \partial_{n}T_\alpha^\alpha=\kappa_\alpha, \quad \partial_{\alpha}T_\alpha^n=-\kappa_\alpha.
\end{equation}
Since $F$ only depends on the eigenvalue of the matrix, $F^{ij}(D^2u)$ and $D^2u$ can be orthogonally
diagonalized simultaneously, and their product is a symmetric matrix. Hence,
\begin{equation}\label{F}
\sum_{j}F^{nj}u_{\alpha j}=\sum_{j}F^{\alpha j}u_{nj}.
\end{equation}
Substitute (\ref{T}) and (\ref{F}) into (\ref{eq:commutator}), we obtain
\begin{equation*}
L(T_\alpha u)=(1+\kappa_\alpha x_n)Lu_\alpha-\kappa_\alpha x_\alpha Lu_n
-\kappa_\alpha g_{p_n}u_\alpha+\kappa_\alpha g_{p_\alpha}u_n.
\end{equation*}
Using (\ref{eq:first-derivative}), we get
\begin{equation*}
|L(T_\alpha u)|\leq C.
\end{equation*}
Similarly, we have $|L(T_\alpha h)|\leq C(1+\tr F)$. Therefore,
\begin{equation}\label{eq:LTw}
\abs{L(T_\alpha w)}\leq C(1+\tr F).
\end{equation}

Consider
\begin{equation*}
Z_\pm=A\widetilde\rho+B\abs x^2\pm T_\alpha w,
\end{equation*}
where $A, B$ are positive constants. By \eqref{eq:Tw-boundary}, $Z_\pm\geq(B-C)|x|^2$ on $\partial\Om\cap B_r(0)$. On $\Om\cap\partial B_r(0)$, using $\widetilde\rho\geq0$ and gradient estimates,
$Z_\pm\geq Br^2-C$. Choose $B$ sufficiently large and we get $Z_\pm\geq0$ on the boundary.
From
\begin{equation*}
L\abs x^2=2\tr F-2g_p\cdot x,
\end{equation*}
we get $\abs{L\abs x^2}\leq C(1+\tr F)$.
Finally choose $A>C(B+1)/c_1$, using \eqref{eq:barrier} and \eqref{eq:LTw},
we have $LZ_\pm\leq0$.
The minimum principle gives $Z_\pm\geq0$ in $\Om\cap B_r(0)$.
$x_0$ lies on the true boundary, where $\widetilde\rho=0$, $T_\alpha w=0$ and
$\abs x^2=0$. Hence $Z_\pm(0)=0$.
Taking the outer normal derivative, we obtain
\begin{equation*}
A\widetilde\rho_\nu(0)\pm(T_\alpha w)_\nu(0)\leq0.
\end{equation*}
Therefore $\abs{(T_\alpha w)_\nu(0)}\leq C$.
At $x_0$, $T_\alpha=\partial_\alpha$ and
$\partial_nT_\alpha^\alpha=\kappa_\alpha$. Thus
\begin{equation*}
(T_\alpha w)_\nu(0)=\kappa_\alpha w_\alpha+(1+\kappa_\alpha x_n)w_{\alpha n}-\kappa_\alpha x_\alpha w_{nn}=w_{\alpha n}(0)+\kappa_\alpha w_\alpha(0).
\end{equation*}
By the boundary condition $w_\alpha(0)=0$, 
\begin{equation*}
\abs{u_{\alpha\nu}(0)}\leq C.
\end{equation*}
Note that the boundary point $x_0$ is arbitrary, so we finally obtain
\begin{equation}\label{eq:mixed-bound}
\abs{u_{\alpha\nu}}\leq C,\quad\alpha<n\quad\text{on }\partial\Om.
\end{equation}

\section{boundary normal $C^2$ estimates}\label{sec:normal}

To deal with the normal derivatives, we use the method in \cite{Trudinger}. For the $(n-1)$-order matrix $T$, denote
\begin{equation*}
T_{\alpha\beta}=(D^2u)_{\alpha\beta},\quad \alpha,\beta\neq n,\quad H(T)=S_{k-1}^\frac{1}{k-1}(T),\quad H^{\alpha\beta}:=\frac{\partial H}{\partial T_{\alpha\beta}}.
\end{equation*}
From $\diag(D^2u,a)\in\Gamma_k^{n+1}$, we obtain
$\diag(T,a)\in\Gamma_{k-1}^n$.
Denote
\begin{equation*}
\widetilde H^{\alpha\beta}(T)=\frac{H^{\alpha\beta}(T)}{\sum_{\gamma=1}^{n-1}H^{\gamma\gamma}(T)}.
\end{equation*}
Then $\widetilde H$ is positive definite with $\tr \widetilde H=1$.
Suppose a curvature matrix $K$ satisfies $K-\delta I\in\Gamma_{k-1}^{n-1}$. Applying
\eqref{eq:support} to $\diag(T,a)$ and $\diag(K-\delta I,0)$ gives
\begin{equation*}
H^{\alpha\beta}(K-\delta I)_{\alpha\beta}\geq\sigma_{k-1}^\frac{1}{k-1}(K-\delta I)>0.
\end{equation*}
Thus,
\begin{equation}\label{eq:dual-curvature}
\widetilde H^{\alpha\beta}K_{\alpha\beta}\geq\delta
\end{equation}
for $K$ at any boundary point.

Choose $x_0$ on the boundary such that
\begin{equation*}
H(T(x_0))=\min_{\partial\Om}H(T).
\end{equation*}
In the neighborhood of $x_0$, take a smooth orthonormal tangential frame
$E_1,\ldots,E_{n-1}$, and represent all tangential matrices in
this frame. We set $K_{\alpha\beta}:=\II(E_\alpha,E_\beta)$. Let $x_0=0$ and $e_n=\nu(0)$.
By concavity, on the boundary nearby we have
\begin{equation*}
0\leq H(T)-H(T(0))\leq H^{\alpha\beta}(T(0))(T_{\alpha\beta}-T_{\alpha\beta}(0)),
\end{equation*}
hence
\begin{equation}\label{eq:frozen-support}
\widetilde H^{\alpha\beta}(T(0))T_{\alpha\beta}\geq \widetilde H^{\alpha\beta}(T(0))T_{\alpha\beta}(0).
\end{equation}
Using $w$ in Section \ref{sec:mixed}, since $w=0$ on the boundary, by \eqref{eq:tangential} we
set
\begin{equation*}
T_{\alpha\beta}=D^2h(E_\alpha,E_\beta)+w_\nu K_{\alpha\beta}.
\end{equation*}
Substituting into \eqref{eq:frozen-support} and then using
\eqref{eq:dual-curvature}, we obtain
\begin{equation}\label{eq:b}
w_\nu\geq b(x):=\frac{\widetilde H^{\alpha\beta}(T(0))T_{\alpha\beta}(0)-\widetilde H^{\alpha\beta}(T(0))D^2h(E_\alpha,E_\beta)}{\widetilde H^{\alpha\beta}(T(0))K_{\alpha\beta}}
\end{equation}
and $w_\nu(x_0)=b(x_0)$.
Note that $D^2h$ and $K_{\alpha\beta}$ are uniformly bounded in $C^2$. With (\ref{eq:dual-curvature})
we have
\begin{equation}\label{eq:b-C2}
\norm b_{C^2}\leq C.
\end{equation}

Next we construct a barrier using derivatives in fixed directions.
Let $\nu_n=\nu\cdot e_n\geq1/2$ on the nearby boundary in sufficiently small neighborhood.
Denote $\pi(x)$ be the nearest-point projection of the tubular neighborhood, and
define
\begin{equation*}
Q(x)=\nu_n(\pi(x))b(\pi(x)).
\end{equation*}
By \eqref{eq:b-C2}, 
\begin{equation}\label{Q}
\norm Q_{C^2}\leq C.
\end{equation}
On $\partial \Omega$ $Dw=w_\nu\nu$, so from (\ref{eq:b})
\begin{equation*}
w_n-Q=\nu_n(w_\nu-b)\geq0,\quad (w_n-Q)(0)=0.
\end{equation*}
Using \eqref{eq:first-derivative}, we get
\begin{equation}\label{eq:L-Psi}
\abs{L(w_n-Q)}=\abs{g_{x_n}+g_zu_n-Lh_n-LQ}\leq C(1+\tr F).
\end{equation}
In the small neighborhood $\Om\cap B_r(0)$, consider
\begin{equation*}
Z=A\widetilde\rho+B\abs{x}^2+(w_n-Q),
\end{equation*}
where $A, B$ are positive constants. On $\partial\Om\cap B_r(0)$ every term is nonnegative. On $\Om\cap\partial B_r(0)$,
by (\ref{eq:C1}) and (\ref{Q}) we choose $B$ sufficiently large so that $Br^2+w_n-Q\geq 0$.
Then choose $A$ sufficiently large, so that \eqref{eq:barrier} and \eqref{eq:L-Psi}
guarantee $LZ\leq0$.
The minimum principle gives $Z\geq0$.
Since $Z(0)=0$, differentiating along the outer normal gives
\begin{equation*}
(w_n-Q)_\nu(0)\leq-A\widetilde\rho_\nu(0)\leq C.
\end{equation*}
Note that $Q$ is defined by projection, so $Q$ is extended as a constant in the normal direction and
$Q_\nu=0$. Hence $(w_n-Q)_\nu(0)=w_{nn}(0)$.
Therefore
\begin{equation*}
u_{\nu\nu}(0)\leq C.
\end{equation*}
On the other hand, $\Delta u>0$ and $\abs{T(0)}\leq C$ give
$u_{\nu\nu}(0)>-\tr T(0)\geq-C$.
Together with (\ref{eq:tangential-bound}) and (\ref{eq:mixed-bound}), we obtain
\begin{equation}\label{d2u0}
\abs{D^2u(0)}\leq C.
\end{equation}

By (\ref{eq:C1}) and $g>0$, there exists $m_0>0$ on the fixed compact set
such that
\begin{equation}\label{gm0}
g(x,u,Du)\geq m_0.
\end{equation}
So with (\ref{d2u0}) we get
\begin{equation*}
S_k(D^2u(0))=g^k(0,u(0),Du(0))\geq m_0^k.
\end{equation*}
Lemma \ref{lem:compact} then guarantees that 
\begin{equation*}
F^{nn}(D^2u(0))\geq \theta>0.
\end{equation*}
Thus there exist $c_0>0$ independent of $u$ and $x_0$ such that
\begin{equation}\label{eq:coefficient-lower}
S_{k-1}(T(x))\geq S_{k-1}(T(0))=kF^{k-1}(D^2u(0))F^{nn}(D^2u(0))\geq c_0,\quad x\in\partial\Om.
\end{equation}

At any boundary point, for $e_n=\nu(x)$, denote
\begin{equation*}
D^2u=\begin{pmatrix}T&\eta_u\\\eta_u^{\mathsf T}&u_{nn}\end{pmatrix},
\quad \eta_u=(u_{1n},\ldots,u_{n-1,n})^{\mathsf T}.
\end{equation*}
An elementary symmetric function is a sum of principal minors of the same
order, and each principal minor contains $u_{nn}$ at most once. Hence
\begin{equation*}
S_k(D^2u)=u_{nn}S_{k-1}(T)+R_a(T,\eta_u).
\end{equation*}
Orthogonally diagonalize $T$ as $\diag(\tau_1,\ldots,\tau_{n-1})$, and write
$\eta_u=(\eta_1,\ldots,\eta_{n-1})^{\mathsf T}$ in the same frame. Then
\begin{equation*}
\begin{split}
R_a(T,\eta_u)={}&\sigma_k(T)+a\sigma_{k-1}(T)\\
&-\sum_{\alpha<n}\eta_\alpha^2
\big[\sigma_{k-2}(\tau|\alpha)+a\sigma_{k-3}(\tau|\alpha)\big].
\end{split}
\end{equation*}
From (\ref{eq:tangential-bound}) and (\ref{eq:mixed-bound}), the tangential block $T$ and the mixed vector $\eta_u$ are already bounded, so $R_a$ is
bounded. $g^k$ is bounded on the fixed compact set. And
\eqref{eq:coefficient-lower} gives a positive lower bound for $S_{k-1}(T(x))$.
Hence
\begin{equation*}
u_{nn}=\frac{g^k-R_a(T,\eta_u)}{S_{k-1}(T)}
\end{equation*}
has a uniform absolute-value bound on the whole boundary.
We obtain
\begin{equation*}
\norm{D^2u}_{L^\infty(\partial\Om)}\leq C.
\end{equation*}

Thus, we obtain the boundary $C^2$ estimate. Combining with the theorem 1.1 (b) in \cite{LiuRen}, we finally prove the global $C^2$ estimate.

\section{the existence and uniqueness}\label{sec:continuity}

Denote
\begin{equation*}
q(x):=F(D^2\underline{u})-g(x,\underline{u},D\underline{u})\geq1.
\end{equation*}
For $0\leq t\leq1$, consider the following path of continuity
\begin{equation}\label{eq:path}
\begin{cases}
F(D^2u)=g_t(x,u,Du):=g(x,u,Du)+(1-t)q(x),&x\in\Om,\\
u=\phi,&x\in\partial\Om.
\end{cases}
\end{equation}
At $t=0$, $\underline{u}$ is the known admissible solution. At $t=1$, (\ref{eq:path}) reduces to the original problem (\ref{eq:problem}).
\par
Since $q$ depends only on $x$,
\begin{equation*}
\begin{gathered}
(g_t)_z=g_z,\quad (g_t)_p=g_p,\quad D_{pp}^2g_t=D_{pp}^2g,\\
D_xg_t=D_xg+(1-t)Dq,\quad
D_{xx}^2g_t=D_{xx}^2g+(1-t)D^2q.
\end{gathered}
\end{equation*}
Therefore, similar to the preceding argument, for any smooth admissible solution $u$ of (\ref{eq:path}), we have
\begin{equation}\label{uc2}
\norm u_{C^2(\bar\Omega)}\leq C,\quad 0\leq t\leq 1.
\end{equation}
\par
From (\ref{gm0}), $g_t^k(x,u,Du)\geq g^k(x,u,Du)\geq m_0^k$. Combining with (\ref{uc2})
and Lemma \ref{lem:compact}, we obtain a fixed compact set $\widetilde K\Subset\widetilde\Gamma_k^n$
containing the ranges of the Hessians of all path solutions, with
\begin{equation*}
\theta_0 I\leq (F^{ij}(N_0))\leq\theta_1 I,\quad N_0\in\widetilde K.
\end{equation*}
In order to apply the uniform elliptic regularity theory defined on the whole
space of matrices, for any $M_0\in\Sym(n)$ define
\begin{equation*}
\widetilde F(M_0)=\inf_{N_0\in\widetilde K}
\big\{F(N_0)+F^{ij}(N_0)(M_0-N_0)_{ij}\big\}.
\end{equation*}
This is the lower envelope of a family of affine functions with bounded
slopes, hence finite, continuous and concave.
By the supporting hyperplane inequality, $\widetilde F=F$ on $\widetilde K$.
For any $P_0\geq0$, the increment of each affine function from $M_0$ to
$M_0+P_0$ lies between $\theta_0\tr P_0$ and $\theta_1\tr P_0$. And the same holds
for the infimum, so $\widetilde F$ is uniformly elliptic on the whole space.
By (\ref{uc2}),
\begin{equation*}
D_{x_i} g_t=g_{x_i}+g_zu_i+g_{p_l}u_{l i}+(1-t)q_i
\end{equation*}
is uniformly bounded on $\OO$.
Therefore $g_t$ has a uniform H\"older bound.
\par
From \cite{SilvestreSirakov}, there exists
$\gamma\in(0,1)$ such that, for each $x_0\in\partial\Omega$, there is a
quadratic polynomial $u_0$ satisfying
\begin{equation*}
\widetilde F(D^2u_0)=g_t(x_0,u(x_0),Du(x_0))
\end{equation*}
and
\begin{equation*}
|u(x)-u_0(x)|\leq C|x-x_0|^{2+\gamma}
\end{equation*}
for $x\in\overline\Omega$ sufficiently close to $x_0$.
Moreover, the map $x_0\mapsto D^2u_0$ has a uniform
$C^{0,\gamma}(\partial\Omega)$ bound.

For an interior point $x$ near the boundary, choose a nearest boundary
point $x_0$ and set $d=|x-x_0|$. Rescale $u-u_0$ on $B_\frac{d}{2}(x)$ and the
corresponding norms of the function and the right-hand side are both $O(d^\gamma)$.
The interior estimate
in~\cite{CaffarelliCabre} therefore yields
\begin{equation*}
|D^2u(x)-D^2u_0|
+d^{\bar\alpha}[D^2u]_{C^{0,\bar\alpha}(B_\frac{d}{4}(x))}
\leq Cd^\gamma
\end{equation*}
for some fixed $0<\bar\alpha\leq\gamma$. Combining this with the
boundary H\"older bound for $D^2u_0$ and the interior estimate
away from the boundary, we obtain
\begin{equation}\label{eq:path-C2alpha}
\|u\|_{C^{2,\bar\alpha}(\overline\Omega)}\leq C.
\end{equation}
Both $\bar\alpha$ and $C$ are independent of $t$ and of the corresponding
solution.

Since $D^2u(\overline\Omega)$
is compactly contained in the admissible branch, the linearization
of \eqref{eq:path} at $u$ is uniformly elliptic.
Interior and boundary Schauder bootstrapping then gives
$u\in C^\infty(\overline\Omega)$.

For $0<\alpha<\bar\alpha$, define
\begin{equation*}
X=\{w\in C^{2,\alpha}(\OO)\big|w|_{\partial\Om}=0\},\quad
Y=C^{0,\alpha}(\OO).
\end{equation*}
Set $u=h+w$ and
\begin{equation*}
\widetilde U=\{w\in X\big|D^2(h+w)(x)\in\widetilde\Gamma_k^n\ \text{for all }x\in\OO\}.
\end{equation*}
$\widetilde U$ is open in $X$.
Denote
\begin{equation*}
G(t,w)=F(D^2(h+w))-g_t(x,h+w,D(h+w))\in Y.
\end{equation*}
\begin{equation*}
\widetilde I=\{t\in[0,1]\big|\eqref{eq:path}\ \text{admits a smooth admissible solution}\}.
\end{equation*}
By the construction of $\underline{u}$, we have $0\in\widetilde I$.

If $t\in\widetilde I$, the Fr\'echet derivative with respect to $w$ at the
corresponding solution is
\begin{equation*}
D_wG(t,w)[\zeta]
=F^{ij}\zeta_{ij}-(g_t)_{p_l}\zeta_l-(g_t)_z\zeta,
\quad\zeta|_{\partial\Om}=0.
\end{equation*}
The operator $D_wG(t,w)$ is a uniformly elliptic linear Dirichlet operator with smooth
coefficients and nonpositive zeroth-order coefficient.
The maximum principle gives injectivity, and classical solvability of the linear Dirichlet problem and the Schauder estimate give surjectivity and a bounded inverse.
It is therefore an isomorphism from $X$ to $Y$.
The implicit function theorem provides a $C^{2,\alpha}$ solution near $t$.
Further improvement of regularity yields a smooth solution,
hence $\widetilde I$ is relatively open in $[0,1]$.

Let $t_j\in\widetilde I$ with $t_j\to t_\infty$, and let $u_j$ be the
corresponding solutions.
By \eqref{eq:path-C2alpha} and the compact embedding, we choose a
subsequence such that
\begin{equation*}
u_j\longrightarrow u_\infty\quad\text{in }C^{2,\alpha}(\OO).
\end{equation*}
The common compact range of the Hessians guarantees that the limit still lies
strictly in $\widetilde\Gamma_k^n$.
Substituting the convergence into the equation and the boundary condition,
we get
\begin{equation*}
F(D^2u_\infty)=g_{t_\infty}(x,u_\infty,Du_\infty),\quad
u_\infty|_{\partial\Om}=\phi.
\end{equation*}
After improving regularity, $u_\infty\in C^\infty(\OO)$, hence
$t_\infty\in\widetilde I$.
Therefore $\widetilde I$ is relatively closed.
Since $[0,1]$ is connected and $\widetilde I$ is nonempty, relatively open and
relatively closed, we have $\widetilde I=[0,1]$.
Then $t=1$ gives the solution required by Theorem \ref{thm:main}.

Let $u_1,u_2$ be two smooth admissible solutions.
Since $f>0$, the Hessians of both lie in $\widetilde\Gamma_k^n$.
Set $w=u_1-u_2$ and $u_s=u_2+sw$ for $0\leq s\leq1$.
By the convexity of $\widetilde\Gamma_k^n$, $D^2u_s\in\widetilde\Gamma_k^n$.
Both solutions satisfy $F(D^2u_i)-g(x,u_i,Du_i)=0$.
Therefore by integrating,
\begin{equation*}
a^{ij}w_{ij}-b^l w_l-cw=0,\quad w|_{\partial\Om}=0,
\end{equation*}
where
\begin{equation*}
a^{ij}=\int_0^1F^{ij}(D^2u_s)\,ds,\quad b^l=\int_0^1g_{p_l}(x,u_s,Du_s)\,ds,\quad
c=\int_0^1g_z(x,u_s,Du_s)\,ds\geq0.
\end{equation*}
On the compact set $\OO\times[0,1]$, $(a^{ij})$ is uniformly positive
definite.
Applying the maximum principle to $w$ and $-w$ separately, we obtain $w=0$.
This completes the proof of Theorem \ref{thm:main}.

\end{document}